\input amssym.def
\input amssym.tex

\input xy
\xyoption{all}

\magnification=\magstephalf
\hfuzz=3pt
\baselineskip=12.8pt
\font\prm=cmr10 at9pt
\font\pit=cmsl10 at9pt
\font\pbo=cmb10 at9pt
\font\psl=cmsl10 at9pt
\font\nrm=cmcsc10 at10pt

\font\bigrm=cmb10 scaled 1200

\def\#{{top}}
\def\stilde{\widetilde}
\def\fini{{$\quad\quad\square$}}
\def\geq{\geqslant}
\def\leq{\leqslant}

\def\ra{\rightarrow}
\def\lra{\longrightarrow}

\def\ext{\hbox{\rm Ext}}
\def\tor{\hbox{\rm Tor}}
\def\hom{\hbox{\rm Hom}}
\def\irr{\hbox{\rm Irr}}

\def\chara{{\rm char}}
\def\proj{{\rm Proj}}
\def\supp{{\rm Supp}}
\def\irr{{\rm Irr}}

\def\R{{\cal R}}

\def\O{{\cal O}}
\def\SI{{{\cal S}_{I}}}
\def\RI{{{\cal R}_{I}}}

\def\reg{{\rm reg}}
\def\ia{{\frak a}}
\def\ib{{\frak b}}

\def\ra{\rightarrow}
\def\lra{\longrightarrow}
\def\ip{{\frak p}}

\def\ig{{\frak g}}
\def\Ip{{\frak P}}
\def\s{{\sigma}}
\def\a{{\alpha}}
\def\b{{\beta}}
\def\g{{\gamma}}
\def\im{{\frak m}}
\def\o{{\omega}}
\def\Sc{{\cal S}}
\def\Oc{{\cal O}}
\def\O{{\Omega}}
\def\cc{{\bf c}}
\def\ol{\overline}
\def\Y{{\cal Y}}

\def\Hc{{H}}
\def\C{{\cal C}}

\ 
\bigskip\bigskip

\centerline{\bigrm LIAISON AND CASTELNUOVO-MUMFORD REGULARITY}
\bigskip\bigskip
\centerline{\nrm By Marc Chardin and Bernd Ulrich\footnote{$^{*}$}{{\prm
Supported in part by the National Science Foundation.\hfill\break
\indent Manuscript received October 27, 2000.\hfill\break}
\indent {\pit American Journal of Mathematics (to appear).}}}
\bigskip
\centerline{\it Dedicated to J{\"u}rgen Herzog on his sixtieth birthday}
\bigskip
\centerline{\vrule height.5pt width5cm}
\bigskip
\bigskip

{\pit Abstract.} {\prm In this article we establish bounds for the 
Castelnuovo-Mumford
regularity of projective schemes in terms of the degrees of
their defining equations. The main new ingredient in our proof is to show
that generic residual intersections of complete intersection
rational singularities again have rational singularities. When applied to
the theory of residual intersections this circle of ideas also sheds new
light on some known classes of free resolutions of residual ideals.}
\bigskip

{\bf Introduction.}
 In this article we prove bounds for the 
Castelnuovo-Mumford regularity of projective schemes in terms 
of the degrees of defining equations, very much in the spirit
of Bertram, Ein, and Lazarsfeld ([BEL]).
Our methods are based on liaison theory. They also lead to
results in positive characteristic and provide information on
defining ideals, even if these are not saturated or unmixed.

The following gives a flavor of one of our main results 
(Theorem 4.7(a)).\smallskip

{\nrm Theorem 0.1.} {\it Let $k$ be a field of characteristic zero 
and ${\Sc}\subset {\bf P}^{n}_{k}$ an equidimensional subscheme of 
codimension $r$ with no embedded components, defined by forms of
degrees $d_{1}\geq \cdots \geq d_{t}\geq 1$. Assume ${\Sc}$ is 
locally a complete intersection except possibly at finitely many
points, and has only rational singularities outside a subscheme
of dimension one. If ${\Sc}$ is not a complete intersection, then
$$
\reg ({\Sc})\leq d_{1}+\cdots +d_{r}-r-1.
$$}

To prove this result we pass to a sufficiently general link ${\Sc}'$ 
of ${\Sc}$ and apply an improved version of Kodaira vanishing to ${\Sc}'$.
This gives a short proof if ${\Sc}$ is smooth of dimension at most three,
since smoothness passes to ${\Sc}'$ ([Hi], [PS]). While this is no
longer true in higher dimensions, we are able to show that a generic 
link of a complete intersection rational singularity again has rational
singularities (Theorem 3.10). We deduce this from a general observation
relating the Rees algebra of an ideal to various residual intersections.
This result also has applications to the theory of residual 
intersections and sheds new light on some known classes of free 
resolutions for such ideals ([BKM], [BE], [KU1]) (Proposition 3.7,
Corollary 3.8, and Remark 3.9).

Our main result in arbitrary characteristic (Theorem 4.7(b))
uses the notion of F-rationality that extends the concept of
rational singularities to any characteristic ([HH1]).
\smallskip 
{\nrm Theorem 0.2.} {\it Let $k$ be any field
and ${\Sc}=\proj (R/I)\subset {\bf P}^{n}_{k}=\proj (R)$ a
locally complete intersection scheme of codimension $r$ with
$1<r<n$. Assume $\Sc$ has at most isolated irrational singularities if
$\chara (k)=0$ or has at most F-rational singularities if $\chara
(k)>0$. Further suppose $I$ is generated by forms of degrees
$d_{1}\geq \cdots \geq d_{t}\geq 1$ and is not a complete
intersection. Then 
$$
\reg (R/I)\leq {{(\dim \Sc +2)!}\over{2}}(d_{1}+\cdots +d_{r}-r-1).
$$}

Notice that $I$ is not required to be saturated or equidimensional. 
To prove the result we induct on the dimension of the scheme by 
passing to the intersection of $\Sc$ with a generic link $\Sc '$.
Two ingredients are used in this reduction: first, a result of K. Smith 
([S2]) that allows one to control the degrees of the equations 
defining $\Sc '$; second, a theorem we prove which says that 
F-rationality passes from $\Sc$ to $\Sc\cap \Sc'$ (Theorem 4.4). 
The latter theorem follows from a broader study of how F-rationality 
behaves under taking generic residual intersections (Theorem 3.13).
\bigskip

{\bf Acknowledgements.} The first author is grateful to Michigan 
State University for its hospitality while this work was completed.
\bigskip

{\bf 1. On the equations defining a link.}
To prove the above theorems we will need to control the degrees of the
equations defining a direct link of a scheme. This is the content of
Theorem 1.7. Its proof uses liaison theory and relies on 
a Kodaira vanishing theorem due to Ohsawa (in characteristic zero) as
well as on work by K. Smith (in positive characteristic). To give a
flavor of the results in this section, we start with the following 
observation:\medskip 

{\nrm Proposition 1.1.} {\it Let $A$ be a standard graded
algebra over a perfect field $k$ such that $A_{\ip}$ is regular for 
every minimal prime $\ip$ of maximal dimension. Let $d:=\dim A$,
$\o :=\o_{A}$ and let $C$ be defined by the sequence
$$
0\ra (\o_{\leq d})\buildrel{\rm can}\over{\lra} \o \lra C \ra 0.
$$
Then ${\rm Supp}(C)\subset {\rm Sing}(A)$.}\medskip

{\it Proof.} Let $\ia \subset A$ be the intersection of all primary
components of $0$ of maximal dimension. Without changing $\o$
or $\dim A$, we may replace $A$ by $A/\ia$ to assume that $A$ is
reduced and equidimensional. Consider the fundamental class
$$
\cc_{A/k}:\wedge^{d}\O_{A/k}\lra \o_{A/k},
$$
a natural $A$-linear map from the module of differential $d$-forms
to the Dedekind complementary module, which is a particular graded
canonical module of $A$ sitting inside the module of meromorphic
$d$-forms ([Elz, p. 34], [KW, 4.11 and 5.13], [Li, 3.1]).
The homomorphism $\cc_{A/k}$ is homogeneous and its cokernel is
supported on the singular locus of $A$ ([KW, 5.13]). Now the assertion
follows since $\wedge^{d}\O_{A/k}$ is generated in degree $d$.\fini
\medskip

Let us recall the definition of the Castelnuovo-Mumford regularity
of a graded module and of a sheaf of modules over an embedded projective 
scheme.\medskip
 
{\it Definition.}  If $S$ is a polynomial ring over a Noetherian ring $B$
and $M$ a finitely generated graded $S$-module, then  
$$
\eqalign{
\reg (M)&:=\min \{ \mu\ |\ H^{i}_{S_{+}}(M)_{>\mu -i}=0,\ \forall i\}\cr
        &=\min \{ \mu\ |\ \tor_{i}^{S}(M,B)_{>\mu +i}=0,\ \forall i\} 
.\cr}
$$

If $\Sc$ is a subscheme of ${\bf P}^{n}_{B}:=\proj (S)$, we set $\reg
(\Sc ):=\reg (S/I)$, where $I$ is the unique saturated ideal such that
$\Sc =\proj (S/I)$. 

If ${\cal M}$ is a coherent sheaf on ${\bf P}^{n}_{B}$ such that
$\Gamma_{\bullet}({\cal M}):=\bigoplus_{\mu \in {\bf Z}}H^{0}({\bf
P}^{n}_{B},{\cal  M}(\mu ))$ is finitely generated over $S$, we 
set $\reg ({\cal M}):=\reg  (\Gamma_{\bullet}({\cal M}))$.
\medskip

{\nrm Corollary 1.2.} {\it Let $R$ be a standard graded Gorenstein 
algebra over a perfect field, and let $\Sc$ be a projective equidimensional
reduced subscheme of ${\cal Z}:=\proj (R)$. Consider a direct link
$\Sc '$ of $\Sc$ in  ${\cal Z}$, given by forms of degrees $d_{1},
\ldots ,d_{r}$, 
and put $\s :=\reg (R)+\sum_{i=1}^{r}(d_{i}-1)$. If $J$ is the defining 
ideal of $\Sc '$, then the scheme defined by $(J_{\leq \s})$ coincides
with $\Sc '$ outside the singular locus of $\Sc$.}
\medskip

{\it Proof.} Let $I$ be the defining ideal of $\Sc$,
$\o :=\o_{R/I}$ and $d:=\dim R/I$. 
We may assume that $\s \geq d_{i}$ for every $i$, since otherwise one
reduces to the case where $r=1$ and $R$ is regular. Let $\ib$ be the
ideal generated by the forms linking $\Sc$ to $\Sc '$. Notice that
$\ib =(\ib_{\leq \s})$.  From the graded isomorphisms $J/\ib\simeq
\hom_{R/\ib}(R/I,R/\ib )\simeq \ext ^{r}_{R}(R/I,R)[-d_{1}\cdots
-d_{r}]\simeq \o [d-\s ]$ we obtain a commutative diagram 
$$
\matrix{
(J/\ib)_{\leq \s} &\buildrel{\sim}\over{\lra}&(\o [d-\s ]_{\leq
\s})&\!\!\!\! =(\o_{\leq d})[d-\s ]\cr 
\downarrow& &\downarrow &\cr
J/\ib &\buildrel{\sim}\over{\lra}&\o [d-\s ]&\cr 
\downarrow& &\downarrow &\cr
D &{\lra}&C[d-\s ],&\cr}
$$
that identifies the cokernels $D$ and $C[d-\s ]$. 

Now by Proposition 1.1, $C$ is
supported on ${\rm Sing}(R/I)$, therefore so is $D$. As $\ib$ is
generated in degree $\leq \s$ the conclusion follows.\fini\medskip

{\it Definition.} A scheme $\Sc$, essentially of finite type over a
field of characteristic zero, has {\it rational singularities} if it
is normal and, for $\Sc '\buildrel{\pi}\over{\lra}\Sc$ a resolution of  
singularities of $\Sc$, one of the following equivalent properties 
holds,\smallskip
{\rm (1)} ${\bf R}^{i}\pi_{*}\Oc_{\Sc'}=0$ for $i>0$,\smallskip
{\rm (2)} $\Sc$ is Cohen-Macaulay and $\pi_{*}\o_{\Sc '}=\o_{\Sc}$.\medskip

Notice that it is implicit in the definition that these properties are
independent of the desingularization. Notice also that $\pi_{*}\o_{\Sc
'}$ is always a subsheaf of $\o_{\Sc}$ and that the normality of $\Sc$
implies that $\pi_{*}\Oc_{\Sc '}=\Oc_{\Sc}$.\medskip

{\it Definition.} A ring $R$ of prime characteristic is
{\it F-rational} if every parameter ideal in $R$ is tightly closed. 
A scheme is F-rational if all its local rings are F-rational.\medskip   

It has not been proved that the F-rational property is stable under
localization in general. However this is known for rings that are
homomorphic images of Cohen-Macaulay rings ([HH1, 4.2(f)]), a
condition always satisfied in our context.

The notion of F-rationality extends to rings essentially of finite 
type over a field of characteristic zero (via reduction modulo $p\gg 0$, 
see  [S1, 4.1] for a precise definition) in which case one talks about
rings of {\it F-rational type}. Deviating slightly from standard  
terminology, we will say that a ring $R$ is of {\it rational type} if
either $R$ is F-rational of prime characteristic, or $R$ is
essentially of finite  type over a field of characteristic zero  and
has F-rational type. 
It was proved by K. Smith in [S1, 4.3] and by Hara and Mehta--Srinivas
in [Har, 1.1] and [MS, 1.1], that in
characteristic zero the notions of rational singularity and
F-rational type coincide. We will denote by $ \irr (\Sc )$ the locus
of the points of non rational type of a scheme $\Sc$. This locus is
closed (at least under an excellence hypothesis satisfied in our
context) by [Ve, 3.5].\medskip

{\nrm Theorem 1.3.} {\it Let $\Sc$ be a projective 
equidimensional scheme over a field of characteristic zero
and ${\cal L}$  an ample invertible sheaf on
$\Sc$, then $H^{i}(\Sc ,\o_{\Sc}\otimes {\cal L})=0$ for $i>\max\{
0,\dim \irr (\Sc )\} $. In particular, if $\dim\irr (\Sc )\leq 0$
then $\reg (\o_{\cal S})=\dim \Sc +1$, for the embedding given by 
any very ample invertible sheaf ${\cal L}$.}\medskip

{\it Proof.} We may suppose that the field is algebraically closed and 
that ${\cal S}$ is reduced. Since the cohomology modules in question 
are epimorphic images of the direct sum of the corresponding cohomology 
modules of the irreducible components of ${\cal S}$, we may further assume 
that ${\cal S}$ is irreducible. 
Let $\Sc '\buildrel\pi\over\lra 
\Sc$ be a desingularization of $\Sc$. One has an exact sequence of
coherent sheaves (see e. g. [Elk, p. 141]), 
$$
(*)\quad 0\ra \pi_{*}\o_{\Sc '}\buildrel{\tau}\over{\lra} \o_{\Sc}\ra
C\ra 0, 
$$
with $\supp (C)\subset \irr (\Sc )$, and from [Oh] 
(or [Ko, 2.1 (iii)]), it follows that 
$$
H^{i}(\Sc , \pi_{*}\o_{\Sc '}\otimes {\cal L})=0
$$
for $i>0$. The exact sequence in cohomology derived from $(*)$ proves
our claim. For the equality  $\reg (\o_{\cal S})=\dim \Sc +1$ also
notice that  $H^{\dim \Sc}(\Sc ,\o_{\Sc})\not= 0$ (see the proof of
Proposition 4.2).\fini\medskip 

{\it Remark 1.4} (i) If $\dim \Sc \geq 1$, ${\cal L}$ is very ample
giving an embedding $\Sc \subset {\bf P}^{n}_{k}$, and $\Sc =\proj
(A)$ with $A=k[X_{0},\ldots ,X_{n}]/I$ for a homogeneous ideal $I$,
then $\reg (\o_{\cal S})=\reg (\o_{A})$. 

(ii) If $\Sc =\proj (A)$, where $A$ is a standard graded
domain of dimension at least two over an algebraically closed field 
of characteristic zero, the 
above proof shows that $\Gamma_{\bullet}(\pi_{*}\o_{\Sc '})$ is a 
submodule of $\o_{A}$ which has  
Castelnuovo-Mumford regularity equal to $\dim A$, and which
coincides with $\o_{A}$ outside the irrational locus of $A$.\medskip

We will need the following partial generalization to positive 
characteristic, which again only needs to be proved for 
reduced and irreducible schemes:\medskip

{\nrm Theorem 1.5} ([S2, 3.2]). {\it Let $\Sc$ be a projective
equidimensional scheme of rational type. If ${\cal L}$  is a globally
generated ample invertible sheaf on $\Sc$, then $\o_{\Sc}\otimes {\cal
L}^{\dim \Sc +1}$ is generated by its global sections.}  
\medskip
{\it Remark 1.6.} 
Even without the rationality assumption, one can prove as in
[S2, the proof of 3.2] using [Ve, 3.9] that $\o_{{\cal S}}\otimes 
{\cal L}^{\dim \Sc +1}$ is generated by its global 
sections on the rational locus of $\Sc$.
Furthermore according to [S3, Theorem 1], if in Theorem 1.5, 
${\cal L}$ is very ample then $\o_{\Sc}\otimes {\cal L}^{\dim \Sc}$ 
is globally generated, unless $\Sc$ is a projective space ${\bf P}$ 
and ${\cal L}= {\cal O}_{\bf P}(1)$. 

\medskip

{\nrm Theorem 1.7.} {\it Let $R$ be a standard graded Gorenstein 
algebra over a field $k$, and let $\Sc$ be a projective equidimensional
reduced subscheme of ${\cal Z}:=\proj (R)$. Consider $\Sc '$, a 
direct link of $\Sc$ in ${\cal Z}$ given by forms of degrees 
$d_{1},\ldots ,d_{r}$ and write $\s :=\reg (R)+\sum_{i=1}^{r} (d_{i}-1)$. 
Denote by $J$ the defining ideal of $\Sc '$ and set $J^{\sharp}:=
(J_{\leq \s})$. 
\smallskip 
{\rm (i)} If $\dim \Sc =1$, then $J=J^{\sharp}$.\smallskip 
{\rm (ii)} If ${\rm char}(k)=0$ and $\dim \irr (\Sc )\leq 0$, 
then $J=J^{\sharp}$.\smallskip 
{\rm (iii)} $J^{\sharp}=J\cap K$ where the scheme defined by $K$ is
supported on the irrational locus of $\Sc$. In particular, if  $\Sc$
is of rational type then $\Sc '=\proj (R/J^{\sharp})$.}\medskip 

{\it Proof.} First notice that in case (i), $\reg (\o_{\Sc})=2$
as $H^{1}(\Sc ,\o_{\Sc}(j))=H^{0}(\Sc ,{\cal O}_{\Sc}(-j))=0$ for
$j>0$. The argument is then the same as in the proof of Corollary 1.2,
replacing Proposition 1.1 by the above equality for (i),
by Theorem 1.3 and Remark 1.4(i) for (ii), and by Remark 1.6 for (iii), 
each time with ${\cal L}=\Oc_{\Sc}(1)$.
\fini\bigskip

{\bf 2. Regularity results for schemes of dimension zero and one.}
If $I$ is an ideal, $I^{\# }$ will denote the unmixed part of $I$,
that is the intersection of the primary components of
height equal to the height of $I$.\medskip

{\nrm Proposition 2.1.} {\it Let $R$ be a standard graded Gorenstein
algebra of dimension $n+1$ over a field and $I$ a homogeneous ideal 
of height $n$.
Assume that $f_{1},\ldots ,f_{n+1}$ are forms of
degrees $d_{1},\ldots, d_{n+1}$ in $I$ so that ${\rm ht} ((f_{1},\ldots
,f_{n+1}):I)\geq n+1$ and
$f_{1},\ldots ,f_{n}$ form a regular sequence. Then either
$$
\reg (R/I)\leq \reg (R)+\sum_{i=1}^{n}(d_{i}-1)
$$
if ${\rm indeg}(I^{\# }/(f_{1},\ldots ,f_{n}))\geq d_{n+1}-1$, or
$
\reg (R/I)\leq \reg (R)+\sum_{i=1}^{n+1}(d_{i}-1)-{\rm
indeg}(I^{\# }/(f_{1},\ldots ,f_{n}))
$
otherwise.}\medskip

{\it Proof.} Write $\ib :=(f_{1},\ldots ,f_{n})$ and $J:=
(f_{1},\ldots ,f_{n+1})$. Notice that $I^{\# }=J^{\# }$ and
$\reg (R/I)\leq \reg (R/J)$. Further, there is an exact
sequence,
$$
0\ra (R/(\ib :J))[-d_{n+1}]
\buildrel{\times f_{n+1}}\over{\lra}
R/\ib \lra R/J\ra 0,
$$
which gives
$$
\reg (R/J)\leq \max \{ \reg (R/\ib ), \reg (R/(\ib
:J))+d_{n+1}-1 \} .
$$

Since $J^{\# }=I^{\# }$ and $\ib :J$ is Cohen-Macaulay,
$$
\eqalign{
\reg (R/(\ib :J))&=\reg (R/\ib )-{\rm indeg}(\ib :(\ib :J)/\ib )
\cr
&=\reg (R/\ib )-{\rm indeg}(I^{\# }/\ib ),\cr}
$$
and this proves our claim.\fini\medskip

{\nrm Proposition 2.2.} {\it Let $R$ be a standard graded Gorenstein
algebra of dimension $n+1$ over a field $k$ and $I$ a homogeneous ideal 
of height $n-1>0$. Assume that $I$ is a complete intersection
locally in codimension $n$ and that $I^{\# }$ is reduced. Let
$f_1,\ldots , f_t$ be forms of degrees
$d_{1}\geq \cdots \geq d_{t}\geq 1$ in $I$ such that 
${\rm ht}((f_1,\ldots , f_t):I)\geq n+1$. Writing $\s :=\reg
(R)+\sum_{i=1}^{n-1}(d_{i}-1)$ one has 
$$
\reg (R/I)\leq \max \{ \s ,3\s -3\},
$$
unless $R=k[X,Y,Z]$ and $I=\ell K$ with $\ell$  a linear form and $K$ 
a complete intersection of 3 forms of degree $d_1-1$ 
(in which case $\reg (R/I)=3\s -2$).}
\medskip
{\it Proof.} We may assume that the ground field is infinite, that
$I=(f_1,\ldots , f_t)$ and that $I$ is not a complete intersection.
There exist forms $\a_{1},\ldots ,\a_{n-1}$ of degrees $d_{1},\ldots
,d_{n-1}$ so that $\a_{1},\ldots ,\a_{n-1}$ form a regular sequence
and $J=(\a_{1},\ldots ,\a_{n-1}):I$ is a geometric link of $I^{\# }$. 
The exact sequence
$$
0\ra R/(\a_{1},\ldots ,\a_{n-1}) \lra R/I \oplus R/J \lra
R/I+J \ra 0
$$
shows that
$$
\reg (R/I)\leq \max \{ \reg (R/(\a_{1},\ldots ,\a_{n-1})),
\reg (R/I+J)\} .
$$

We first assume that $R$ is not regular or that $n-1\geq 2$. In this
case $\s \geq d_{n}$. By Theorem 1.7(i), $J=(J_{\leq \s})$. Let $\Ip$
be the finite set of  primes of
height $n$ containing $I+J$. For $\ip \in \Ip$ either $I_{\ip}$ is
a complete intersection of height $n-1$ and $(I+J/I)_{\ip}$ 
is cyclic, or $I_{\ip}$ is a complete intersection of height $n$  and
$(I+J/I)_{\ip}=0$. Since furthermore $\s\geq d_{n}$, we can choose forms
$\b_{1},\ldots ,\b_{n-1}$ in $I$ of degrees $d_{1},\ldots ,d_{n-1}$ and
a form $\b_{n}$ in $I+J$ of degree $\s$ so that  $\b_{1},\ldots ,\b_{n}$
form a regular sequence, $(I+J)/(\b_{1},\ldots ,\b_{n})$ is generated
in degrees $\leq\s$, and $(I+J)_{\ip}=(\b_{1},\ldots
,\b_{n})_{\ip}$ for every $\ip \in \Ip$. Now there exists a form
$\b_{n+1}$ of degree $\s$ in $I+J$ such that ${\rm ht} ((\b_{1},\ldots
,\b_{n+1}):(I+J))\geq n+1$. Therefore
$$
\reg (R/I+J)\leq \max\{ \s +\s -1, \s +2(\s -1)-1\} 
$$
by Proposition 2.1. Notice also that
$2\s -1\leq  \max \{ \s , 3\s -3\}$ for a non negative integer $\s$.

Next, assume that $R$ is regular and $n-1=1$, hence $R=k[X,Y,Z]$. 
In this case $I= \ell \, K$ with $\ell$ a form of degree $s \geq 1$
and $K$ a homogeneous ideal of height 2 or 3. If $K$ is not a complete 
intersection of height 3, then $\reg (R/K) \leq \max \{
d_1-s-1+d_2-s-1,d_1-s-1 +d_2-s-1+d_3-s-1-1\}$. 
This follows from Proposition 2.1 if $K$ has height 2 and is obvious 
otherwise. Thus
$\reg (R/I)=s+\reg (R/K) \leq \max\{ 2\s -1, 3\s -3 \}$ as asserted. 
If on the other
hand $K$ is a complete intersection of 3 forms of degrees
$s_i \leq d_i -s$, we obtain
$\reg (R/I)= s+s_1-1+s_2-1+s_3-1$, which is $\leq 3\s -3$ unless $s=1$
and $s_1=s_2=s_3=d_1 -1$.
\fini\bigskip

{\bf 3. Controlling residual intersections via blowingup or deformation.}
This section contains the main technical results of the paper. We 
investigate the behavior of rational and F-rational singularities
under generic residual intersections. Two methods are applied, the study
of residual intersections via projecting free resolutions of blow-up
algebras, and the use of deformation theory to reduce to the ``generic'' 
case. The first approach does not yield results about F-rationality, but 
gives insight into resolutions of residual intersections of ideals that 
are not necessarily complete intersections. The second method on the
other hand, allows us to treat F-rationality, but only in the context of
residual intersections of complete intersections.\medskip

We start with a very general lemma from homological algebra.\medskip

{\nrm Lemma 3.1.} {\it Let $R$ be a commutative ring, 
$S:=R[T_{1},\ldots ,T_{r}]$,
and $(D_{\bullet},d_{\bullet})$ a graded complex of $S$-modules with 
$D_{i}=0$ for $i<0$ 
and $D_{i}=\bigoplus_{j\in {\bf Z}}S[-j]^{\b_{ij}}$. Set $b_{i}:={\rm sup}
\{j\; |\;  \b_{ij}\not= 0\}$ and $\ig :=(T_{1},\ldots ,T_{r})$. 
Assume that $H_{q}(({D_{\bullet}})_{T_{\ell}})=0$ for $q>0$ and all
$\ell$, and set  $M:=H_{0}(D_{\bullet})$. Then 
\medskip
{\rm (i)} $H_{p}(D_{\bullet})\simeq H_{r+p}(\Hc^{r}_{\ig}
(D_{\bullet}))$ for $p>0$, 
and $\Hc^{q}_{\ig}(M)\simeq  H_{r-q}(\Hc^{r}_{\ig}(D_{\bullet}))$ 
for $q\geq 0$,
\medskip
{\rm (ii)} $H_{p}(D_{\bullet})_{\nu}=0$ for $\nu >b_{r+p}-r$ and $p>0$,
\medskip
{\rm (iii)} $\Hc^{q}_{\ig}(M)_{\nu}=0$ for $\nu >b_{r-q}-r$ and $q\geq 0$.
\medskip
{\rm (iv)} Let $k\geq 0$ and $t\geq 0$. If $\b_{ij}=0$ for $i<k+r$ and
$j\geq k+r$ and for $k<i< k+\max\{ 2,t\}$ and $j\leq k$, one has an
exact sequence,
$$
\Hc^{r}_{\ig}(D_{k+r+t})_{k}\ra\cdots \ra \Hc^{r}_{\ig}(D_{k+r})_{k}
\buildrel{\tau}\over{\lra}
(D_{k})_{k}\ra \cdots \ra (D_{0})_{k}\buildrel{can}\over{\lra}
 (M/\Hc^{0}_{\ig}(M))_{k}\ra 0.
$$
}

{\it Proof.} Consider the two spectral sequences arising from the
double complex ${\cal C}^{\bullet}_{\ig}(D_{\bullet})$, where 
${\cal C}^{\bullet}_{\ig}(N)$ denotes the \v Cech complex on $N$ relatively
to the ideal $\ig$. Note that as $\ig$ is generated by a regular sequence,
the \v Cech cohomology is the same as the local cohomology. From the 
hypotheses one immediately gets,\smallskip

$
'{\rm (1)}\quad
'E_{1}^{pq}=\left\{ \eqalign{H_{p}(D_{\bullet})&\ \hbox{\rm if}\ q=0,\cr 
                        \C ^{q}_{\ig}(M)&\ \hbox{\rm if}\ p=0,\cr 
                              0\quad&\ \hbox{\rm else},\cr}\right. $\smallskip

$
'{\rm (2)}\quad
'E_{2}^{pq}={'E}_{\infty}^{pq}=\left\{ \eqalign{H_{p}(D_{\bullet})&\ \hbox{\rm
if}\ q=0\ \hbox{\rm and}\ p\not=0,\cr 
                        \Hc _{\ig}^{q}(M)&\ \hbox{\rm if}\ p=0,\cr 
                              0\quad&\ \hbox{\rm else},\cr}\right. $\smallskip

$''{\rm (1)}\quad ''E_{1}^{pq}=\left\{ \eqalign{\Hc
^{r}_{\ig}(D_{p})&\ 
\hbox{\rm if}\ q=r,\cr 
                              0\quad&\ \hbox{\rm else},\cr}\right.$\smallskip

$''{\rm (2)}\quad ''E_{2}^{pq}={''E}_{\infty}^{pq}=\left\{ 
\eqalign{H_{p}(\Hc ^{r}_{\ig} 
(D_{\bullet}))&\ \hbox{\rm if}\ q=r,\cr 
                              0\quad&\ \hbox{\rm else}.\cr}\right.$\smallskip

This immediately gives (i). Also (ii) and (iii) follow as the $''E_{1}^{r+p,r}$ 
and $''E_{1}^{r-q,r}$ terms are zero in the considered degrees.

For (iv), notice that our hypotheses imply that
$(''E_{1}^{ir})_{k}=0$ for $i<k+r$ and $(D_{k+i})_{k}=0$ for 
$0<i<\max\{ 2,t\}$. This shows that the given 
complex is exact to the left of $\tau$, and between $\tau$ and the map
$can$. It also implies that 
$H_{k+r}(\Hc ^{r}_{\ig}(D_{\bullet}))_{k}={\rm coker}(\Hc ^{r}_{\ig}(
d_{k+r+1}))_{k}$. Furthermore $({'E}_{2}^{k0})_{k}=H_{k}(D_{\bullet})_{k}
={\rm ker}(d_{k})_{k}$ and $\Hc ^{0}_{\ig}(M)_{k}=0$ if $k>0$, whereas 
$({'E}_{2}^{00})_{0}=\Hc ^{0}_{\ig}(M)_{0}$ and $(D_{0})_{0}=M_{0}$ if $k=0$.
Therefore, defining $\tau$ as composition of canonical maps
$$
\Hc ^{r}_{\ig}(D_{k+r})_{k}
\lra
H_{k+r}(\Hc ^{r}_{\ig}(D_{\bullet}))_{k}
\buildrel{\delta}\over{\lra}
('E_{2}^{k0})_{k}
\hookrightarrow 
(D_{k})_{k}
$$
proves our claim, as $\delta$ is an isomorphism.\fini
\medskip

In our situation, we will be interested in the more specific result 
below, that essentially follows from the one above. For a ring $R$, an
$R$-module $M$  and a sequence $\g$ of elements of $R$, we will denote
in the sequel by $K_{\bullet}(\g \ ;M)$ the Koszul complex of
$\g$ on $M$ and by $H_{i}(\g\, ;M)$ the $i$-th homology module of 
this complex. \medskip

{\nrm Proposition 3.2.} {\it  Let $R$ be a commutative ring, 
$S:=R[T_{1},\ldots ,T_{r}]$,
$\ig :=(T_{1},\ldots ,T_{r})$, and $M$ a graded $S$-module.
Assume that $M$ has a free $S$-resolution 
$F_{\bullet}$ such that $F_{i}=\bigoplus_{j=0}^{i}S[-j]^{\b_{ij}}$.
Let $\g =\g_{1},\ldots ,\g_{s}$ be a sequence of elements in $S_{1}$ 
such that $H_{i}(\g\, ;M_{T_{\ell}})=0$ for $i>0$ and all $\ell$, let
$J$ be the ideal they generate in $S$ and $D_{\bullet}:=F_{\bullet}
\otimes_{S}K_{\bullet}(\g\, ;S)$. Then
\medskip

{\rm (i)} $H_{i}(\g\, ;M)\simeq H_{r+i}(\Hc^{r}_{\ig}(D_{\bullet}))$ 
for $i>0$, and $\Hc^{i}_{\ig}(M/JM)\simeq
H_{r-i}(\Hc^{r}_{\ig}(D_{\bullet}))$ for $i\geq 0$, 
\medskip

{\rm (ii)} $\Hc^{i}_{\ig}(M/JM)$ is concentrated in degrees $\leq -i$
for all $i$,
\medskip

{\rm (iii)}  $H_{i}(\g\, ;M)$ is concentrated in degree $i$ for $i>0$,
\medskip

{\rm (iv)} $\Hc^{0}_{\ig}(M/JM)_{0}=(0:_{M/JM}\ig )_{0}=
J_{1}M_{0}:_{M_{0}}\ig_{1}$, 
\medskip
{\rm (v)} 
$$
\Gamma (k):=\Gamma ({\rm Proj}(S),\stilde{M/JM}(k))=\left\{ 
\eqalign{
M_{0}/(J_{1}M_{0}:_{M_{0}}\ig_{1})&\ \hbox{\rm if}\ k=0,\cr
         (M/JM)_{k}\quad\quad &\ \hbox{\rm if}\ k>0.\cr}\right.
$$\smallskip

{\rm (vi)} Let $k\geq 0$. If $\b_{ij}=0$ for $i\not= j\leq k$, 
then one has an acyclic complex of free 
$R$-modules,
$$
\cdots \ra \Hc^{r}_{\ig}(D_{k+r+1})_{k}\ra\Hc^{r}_{\ig}(D_{k+r})_{k}
\buildrel{\tau}\over{\lra}
(D_{k})_{k}\ra \cdots \ra (D_{0})_{k}\ra 0
$$
that resolves $\Gamma (k)$.
\medskip

{\rm (vii)}  Let $k\geq 0$. There exists a morphism of complexes 
$u: \Hc^{r}_{\ig}(D_{\bullet})_{k}\{ r\} \lra (D_{\bullet})_{k}$ 
such that $H_{i}(u)$ is an isomorphism for $i>0$ and a monomorphism
for $i=0$, $H_{i}(u)=0$ for $i\not= k$, and $H_{k}(u)$ is the canonical  
map $H_{k+r}(H_{\ig}^{r}(D_{\bullet})_{k})
\lra H_{k}((D_{\bullet})_{k})$. Therefore the mapping 
cone of $u$ is a free $R$-resolution of $\Gamma (k)$.
}
\medskip

{\it Proof.} For (i) and (ii), notice that $D_{\bullet}$ satisfies 
the hypotheses of Lemma 3.1, and, as $F_{\bullet}$ resolves $M$, $H_{i}
(D_{\bullet})\simeq H_{i}(\g\, ;M)$. Lemma 3.1 also shows that $H_{i}
(\g\, ;M)_{\nu}=0$ for $\nu >i$ if $i\not= 0$, and as $M$ is generated 
in degree $0$ and the $\g_{i}$'s are of degree $1$, $K_{i}(\g\, ;M)_{\nu}=0$
for $\nu <i$. Now (iii) follows. 

To prove (iv) note that one always has $0:_{M/JM}\ig \subset 
\Hc^{0}_{\ig}(M/JM)$ and that $\ig [\Hc^{0}_{\ig}(M/JM)]_{0}\subset 
[\Hc^{0}_{\ig}(M/JM)]_{\geq 1}=0$ by (ii), thereby proving the other 
inclusion. Part (v) follows directly from (ii) and (iv),
and part (vi) is a direct consequence of Lemma 3.1(iv). 

Concerning (vii), 
set $L_{\bullet}:=\Hc^{r}_{\ig}(D_{\bullet})_{k}\{ r\}$ and $P_{\bullet}:=
(D_{\bullet})_{k}$. By (i) and (iii), $L_{\bullet}$ and
$P_{\bullet}$ have no homology in homological degrees $>k$ and 
there is a canonical map $H_{k}(L_{\bullet})\buildrel
{\delta}\over{\ra}H_{k}(P_{\bullet})$ that is an isomorphism for $k>0$
and an embedding for $k=0$. Now the truncated complexes $\cdots \ra
L_{k+1}\buildrel{\phi}\over{\ra}L_{k}\ra 0$ and  $\cdots \ra
P_{k+1}\buildrel{\psi}\over{\ra} P_{k}\ra 0$ are free $R$-resolutions
of ${\rm Coker}(\phi )$ and ${\rm Coker}(\psi )$ respectively. 

As $L_{k-1}=0$, the canonical monomorphism $H_{k}(L_{\bullet})
\buildrel{i}\over{\lra}
{\rm Coker}(\phi )$ is an isomorphism. Now the map 
$$
\overline u :{\rm Coker}(\phi )\buildrel{i^{-1}}\over{\lra}H_{k}(L_{\bullet})
\buildrel{\delta}\over{\ra} H_{k}(P_{\bullet})
\buildrel{can}\over{\lra} {\rm Coker}(\psi )
$$ 
can be lifted to a 
morphism between the truncated complexes, which, extended by $0$ to 
a morphism from $L_{\bullet}$ to $P_{\bullet}$, satisfies our requirements.
Notice that $L_{\bullet}$ is trivial in homological degrees $<k$ and 
that by (iii), $P_{\bullet}$ has only homology in homological degrees 
$0$ and $k$. The zeroth homology of the mapping cone of $u$ is given
by (v) if $k>0$ and by (i) and (v) if $k=0$.\fini
\medskip

{\it Remark 3.3.} Assume $R$ is an $\ell$-algebra, graded by a finitely 
generated free abelian group $G$, each $R_{i}$ is a finite free 
$\ell$-module, $M$ is bigraded and the elements $\g_{i}$ are
bi-homogeneous of degree, say, $(e_{i},1)$. Then 
the complexes that appear in Proposition 3.2 are also bigraded. Writing 
$\chi_{S[-k]}(X,Y):=Y^{k}\sum_{i\in G,j\in {\bf Z}}(\hbox{\rm rank}_{\ell}R_{i})
X^{i}Y^{j}$, and 
$\chi_{C_{\bullet}}:=\sum_{p}(-1)^{p}\chi (C_{p})$
for any bigraded complex $C_{\bullet}$ of free $S$-modules, one has
$$
\chi_{D_{\bullet}}(X,Y)=\chi_{F_{\bullet}}(X,Y)\prod_{i=1}^{s}
(1-X^{e_{i}}Y).
$$
From Proposition 3.2(vii) one sees that the Hilbert-Poincar{\'e} series
of $\Gamma (\hbox{---})$ is 
$$
\chi_{D_{\bullet}}(X,Y)-(-Y)^{-r}\chi_{D_{\bullet}}
(X,Y^{-1}),
$$ 
which in turn can be expressed in terms of 
$\chi_{F_{\bullet}}(X,Y)$.\medskip

The following proposition will be a key to passing rational
singularities from the blow-up to a generic residual
intersection.\medskip 

{\nrm Proposition 3.4.} {\it Let $R$ be a local ring essentially of finite 
type over a field of characteristic zero, and let $A$ be a standard graded 
$R$-algebra
with Castelnuovo-Mumford regularity $0$ such that ${\rm Proj}(A)$
has rational singularities. If $x_{1},\ldots ,x_{r}$ generate 
$A_{1}$ set $g _{i}:=\sum_{j=1}^{r}U_{ij}x_{j}\in A_{1}\otimes_{R}
R'$, where $R':=R[U_{ij}]$ is a polynomial ring and $1\leq i\leq s$. 
Then the ring 
$$
R'/((g _{1},\ldots ,g _{s})R':_{R'}A_{1})
$$
has rational singularities.}
\medskip

{\it Proof.} Let $A':=A\otimes_{R}R'$, $S':=R'[T_{1},\ldots ,T_{r}]$
and $\ig := (T_{1},\ldots ,T_{r})$. 
The natural homogeneous map $S'\ra A'$ sending $T_{i}$ to $x_{i}$ makes 
$A'$ a graded $S'$-module. The 
forms $\g_{i}:=\sum_{j=1}^{r}U_{ij}T_{j}\in S_{1}'$ satisfy the hypotheses of
Proposition 3.2. By assumption $A'$ has an $S'$-resolution $F_{\bullet}$ 
satisfying the hypotheses of the same proposition. We set 
$J:=(\g _{1},\ldots ,\g _{s})S'$ and notice that $JA'=(g _{1},\ldots
,g _{s})A'$.  

Note that ${\cal Z}:={\rm Proj}(A'/JA')$ has a natural structure
of vector bundle over ${\rm Proj}(A)$ so that ${\cal Z}$ has rational 
singularities as well. Now consider the projection $\pi : {\cal Z}\ra 
{\rm Spec} (R')$. One has\smallskip

(1) ${\bf R}^{p}\pi_{*}{\cal O}_{\cal
Z}=(H^{p+1}_{\ig}(A'/JA')_{0})^{\sim}=0$  
for $p>0$, by Proposition 3.2(ii), \smallskip

(2) $\pi_{*}{\cal O}_{\cal Z}=(R'/((g _{1},\ldots ,g
_{s})R':_{R'}A_{1}))^{\sim}$ by  Proposition 3.2(v). 

By the Leray spectral sequence ${\bf R}^{p}\pi_{*}({\bf R}^{q}\xi_{*}
\Oc_{\stilde{\cal Z}})\Rightarrow {\bf R}^{p+q}(\pi\circ\xi
)_{*}\Oc_{\stilde{\cal Z}}$ where $\stilde{\cal
Z}\buildrel{\xi}\over{\lra}{\cal Z}$ is a resolution of singularities
of ${\cal Z}$, we get the result.\fini 
\medskip

{\it Remark 3.5.} In any characteristic, one can define a scheme $\Sc$
to have rational singularities if there exists a desingularization
${\Sc '}\buildrel{\pi}\over{\lra} \Sc$ such that ${\bf
R}^{\bullet}\pi_{*}{\cal O}_{{\Sc '}}={\cal O}_{\Sc}$ and ${\bf
R}^{\bullet}\pi_{*}\o_{{\Sc '}}=\o_{\Sc}$. Using this definition one
can show, along the same lines as in the above proof, that Proposition
3.4 holds without the hypothesis that $R$ is essentially of finite
type over a field of characteristic zero. \medskip

In the sequel, if $M$ is a finitely generated module, $\mu (M)$ will 
denote the minimal number of generators of $M$. Let us recall some 
definitions and properties related to blow-up algebras.
\medskip
{\it Definition.} Let $R$ be a Noetherian local ring,  
$z :=z_{1}, \ldots ,z_{s}$ a sequence of elements of $R$,
and $J$ the $R$-ideal they generate.

(i) $z$ is a {\it $d$-sequence} if $\mu (J)=s$ and $((z_{1},\ldots ,z_{j}):
z_{j+1})\cap J=(z_{1},\ldots ,z_{j})$ for $0\leq j<s$,\smallskip

(ii) $z$ is a {\it proper sequence} if $z_{j+1}H_{i}({\rm K}_{\bullet}
(z_{1},\ldots ,z_{j};R))=0$ for $0\leq j<s$ and $i>0$, \smallskip

(iii) $J$ has {\it sliding depth} if ${\rm depth}H_{i}(z ;R)\geq \dim
R-s+i$ for all $i$, \smallskip 

(iv) $J:I$ is an {\it $s$-residual intersection} of the $R$-ideal $I$
properly containing $J$ if ${\rm ht}(J:I)\geq s$. The residual
intersection is called {\it geometric} if in addition  ${\rm
ht}(I+(J:I))> s$. 
\medskip

{\nrm Proposition 3.6.} {\it Let $R$ be a Noetherian 
local ring, $I\subset R$ a proper ideal and $z_{1},\ldots ,z_{r}$  
minimal generators of $I$. Set $S:=R[T_{1},\ldots ,T_{r}]$,
$\SI :={\rm Sym}_{R}(I)$, 
and $\R _{I}:=R\oplus I\oplus I^{2}\oplus \cdots$.
The generators of $I$ define a natural commutative diagram of homogeneous
epimorphisms of degree 0,
$$
\xymatrix{ &S\ar[dl]_{\phi}\ar[dr]^{\psi}& \\
{\SI}\ar[rr]^{\a}& &\R _{I}\\}
$$
and\smallskip 

{\rm (i)} $\a$ is an isomorphism  if and only if $\a\otimes_{R}R/I$ is,\smallskip  

{\rm (ii)} $\ker (\a )$ is nilpotent if and only
$\ker (\a \otimes_{R}R/I)$ is,\smallskip  

{\rm (iii)} If $z_{1},\ldots ,z_{r}$ is a proper sequence then the 
Castelnuovo-Mumford regularity of ${\SI}$ is $0$,
and the converse holds if $R$ has infinite residue field and
$z_{1},\ldots ,z_{r}$ are general,\smallskip   

{\rm (iv)} If $z_{1},\ldots ,z_{r}$ is a $d$-sequence then the 
Castelnuovo-Mumford regularity of $\R _{I}\otimes_{R}R/I$ is $0$,
and the converse holds if $R$ has infinite residue field and
$z_{1},\ldots ,z_{r}$ are general,\smallskip   

{\rm (v)} If $z_{1},\ldots ,z_{r}$ is a $d$-sequence then $\a$ is an
isomorphism, 
\smallskip 

{\rm (vi)} $z_{1},\ldots ,z_{r}$ is a regular sequence if and only if
$\psi\otimes_{R}R/I$ is an isomorphism.}
\medskip

For part (i) see [HSV, 3.1], for (ii) [HSV, 3.2], for (iii) and (iv)
[HSV, 12.7, 12.8, 12.10], and for (v) [Hu, 3.1] or [Va, 3.15]. 

\medskip

{\nrm Proposition 3.7.} {\it Let $R$ be a Noetherian local ring,
$I\subset R$ an ideal generated by a proper sequence of length $r$,
and $\SI :={\rm Sym}_{R}(I)$. For $S:=R[T_{1},\ldots ,T_{r}]$ the map 
$S\ra \SI$ sending $T_{i}$ to the $i$-th generator of $I$ makes $\SI$ 
a graded $S$-module. 

Let $J\subset I$ be an ideal generated by a sequence 
$\g :=\g_{1},\ldots ,\g_{s}$. 
\smallskip

{\rm (a)} Assume ${\rm ht}I>0$, $R$ is universally 
catenary, and $\proj ({\SI})$ is Cohen-Macaulay and equidimensional.
Then $\g \subset ({\SI})_{1}$ form a regular sequence locally on 
$\proj ({\SI})$ if and only if ${\mu}((I/J)_{\ip})\leq 
\dim R_{\ip}-s+1$ for every $\ip\in V(J:I)$. In this case 
${\rm ht}(J:I)\geq s$.
\smallskip
{\rm (b)} Assume $\g \subset ({\SI})_{1}$ form a regular sequence
locally on $\proj ({\SI})$.
\smallskip

{\rm (i)} ${\rm pd}_{R}\, R/(J:I)\leq {\rm pd}_{S}\, {\SI}+s-r+1$. If $R$ 
is regular, then 
$$
{\rm depth}\, R/(J:I)\geq {\rm depth}
\, {\SI}-s-1.
$$
\indent {\rm (ii)} If ${\rm ht}I>0$, $J\not= I$, $R$ is Cohen-Macaulay, 
and $\SI$ is a
perfect $S$-module, then $J:I$ is a perfect ideal of grade $s$.
\smallskip
{\rm (iii)} Let $k>0$ and $t\geq 0$. Assume that ${\rm pd}_{R}\, 
{\rm S}_{j}(I)\leq t+j$ for $1\leq j\leq k$. Then,
${\rm pd}_{R}\, {\rm S}_{k}(I/J)\leq \max \{
t+k,{\rm pd}_{S}\, {\SI}+s-r+1\}$. In particular, if $R$ is regular, then
$$
{\rm depth}\, {\rm S}_{k}(I/J)\geq \min \{\dim R-t-k, {\rm depth}\, 
{\SI}-s-1\} .
$$}
\indent {\it Proof.} (a) First notice that if ${\rm ht}I>0$ and 
$\proj (\SI )$ 
is equidimensional, then $\SI$, $\RI$, and hence $R$ are equidimensional,
and $\dim \SI =\dim \RI =\dim R+1$. Since $\proj (\SI )$ is Cohen-Macaulay,
$\g$ form a regular sequence on   $\proj (\SI )$ if and
only if ${\rm ht}((\g )\SI :(\SI_{+})^{\infty})\geq s$ (here and in what
follows, $(\g )$ denotes the ideal generated by the $\g_{i}$'s as elements
of $(\SI )_{1}$).
One has $(J:_{R}I,(\g ))\SI \subset (\g )\SI :(\SI_{+})^{\infty}$, and
by Proposition 3.2 (ii), (iv) and Proposition 3.6(iii) equality holds
if $\g$ is a regular sequence on  $\proj (\SI )$. Therefore the above
height condition is equivalent to ${\rm ht}(J:_{R}I,(\g ))\SI \geq
s$. As $\SI$ is an equidimensional and catenary positively graded ring
over a local ring, the last inequality means that $\dim
\SI/(J:_{R}I,(\g ))\leq \dim \SI -s=\dim R-s+1$. But $\SI
/(J:_{R}I,(\g ))\simeq {\rm Sym}_{R/J:I}(I/J)$ and by [HR, 2.6],  
$$
\dim {\rm Sym}_{R/J:I}(I/J)=\max\{ \dim R/\ip +\mu ((I/J)_{\ip})\ 
|\ \ip\in V(J:I)\}.
$$
Now the asserted equivalence follows because $R$ is an equidimensional and 
catenary local ring. Furthermore if either condition holds, then for every
$\ip\in V(I:J)$, $ \dim R_{\ip} \geq s+\mu ((I/J)_{\ip})-1\geq s$.

(b) Proposition 3.6(iii) shows that Proposition 3.2 applies. 
By Proposition 3.2(v), parts (vi) and (vii) of that proposition imply
the present assertions (i) and (iii), 
respectively. For (ii), notice that if ${\rm ht}I>0$, and $R$ and $\SI$ are
Cohen-Macaulay, then $\dim \SI =\dim R+1$ and ${\rm grade}(J:I)\geq s$ by 
(a). \fini\medskip

{\nrm Corollary 3.8.} {\it Let $R$ be a regular local ring and $I\subset R$ 
an ideal satisfying sliding depth such that $\mu (I_{\ip})
\leq \dim R_{\ip}+1$ for every $\ip \in V(I)$. Then any 
geometric $s$-residual intersection of $I$ is a perfect ideal of 
grade $s$.}\medskip

{\it Proof.} We may assume that $R$ has infinite residue field
and ${\rm ht}I>0$. As $I$ satisfies sliding depth and $\mu (I_{\ip})
\leq \dim R_{\ip}+1$ for every $\ip \in V(I)$, $\SI$ is Cohen-Macaulay and 
$I$ is generated by a proper sequence 
$z_{1},\ldots ,z_{r}$ ([HSV, 5.3, 10.1, 12.9]). Let $J:I$ be a 
geometric $s$-residual intersection. Write $J=(\g_{1},\ldots ,\g_{s})$
and $\g_{i}:=\sum_{j=1}^{r}a_{ij}z_{j}$ with $a_{ij}\in R$. Set 
$g_{i}:=\sum_{j=1}^{r}U_{ij}z_{j}$ where $U_{ij}$ are variables, and
let $\im$ be the maximal ideal of $R$. Write $R':=R[U_{ij}]_{(\im ,U_{ij}
-a_{ij})}$, $J':=(g_{1},\ldots ,g_{s})R'$, and let $\pi :R'\ra R$ be the
specialization map sending $U_{ij}$ to $a_{ij}$. Notice that $\pi (J')=J$. 
The generic elements $g_{1},\ldots ,g_{s}$ form a regular sequence on
$\proj ({\cal S}_{IR'})$. Hence by Proposition 3.7(b)(ii), $J':IR'$ is
a perfect ideal of grade $s$. 
Since $J:I=\pi (J'):\pi (IR')$ is a geometric residual intersection, 
[HU3, 4.2(ii)] shows that $J:I$ is perfect of grade $s$ as well.
\fini\medskip

{\it Remark 3.9.} (i) In the setting of Proposition 3.7(b) and Corollary
3.8, the complexes of Proposition 3.2 (vi) and (vii) provide free resolutions
of $R/(J:I)$ and of ${\rm S}_{k}(I/J)$, $k>0$, respectively.

(ii) If $I$ is generated by a $d$-sequence, one can apply Proposition 3.2 
to ${\cal G}_{I}:={\cal R}_{I}\otimes_{R}R/I$ instead of $\SI$. This gives 
results for $R/((J+I^{2}):I)$ and $I^{k}/(JI^{k-1}+I^{k+1})$, $k>0$, 
that are similar to the ones of Proposition 3.7.
 
(iii) Let $X$ be a Cohen-Macaulay scheme, let $S\subset Y\subset X$ be
closed subschemes with ${\rm codim}_XS>0$,  
and let $S'\subset X$ be the subscheme with ${\cal I}_{S'}=
{\cal I}_{Y}:{\cal I}_{S}$. We can compare $S'$ to the residual 
scheme proposed in [Fu, 9.2.2]. Indeed, let $\pi :\widetilde{X}:=
{\rm Bl}_{S}X\rightarrow X$, $E:=\pi^{-1}S$ and 
${\cal Z}\subset \widetilde{X}$ the subscheme with ${\cal I}_{{\cal
Z}}={\cal I}_{\pi^{-1}Y}{\cal I}_{E}^{-1}$. If ${\cal I}_{S}$
satisfies sliding depth, $\mu ({\cal I}_{S,x})\leq \dim {\cal
O}_{X,x}$ for every $x\in S$, $\mu ({\cal I}_{Y,x})\leq s$ for every
$x\in Y$, and $\mu ({\cal I}_{S,x}/{\cal I}_{Y,x})\leq \dim {\cal
O}_{X,x}-s+1$ for every $x\in S'$, then Propositions 3.2(v), 3.6, 
3.7(a) and [HSV, 5.3, 9.1, 12.9] show that 
${\cal O}_{S'}=\pi_{*}{\cal O}_{\cal Z}$.\medskip

{\nrm Theorem 3.10.} {\it Let $R$ be a local ring essentially of 
finite type over a field of characteristic zero, and $I=(z_{1},\ldots
,z_{t})\subset R$  an ideal generated by a proper sequence.  
For $1\leq i\leq s$, let $g _{i}:=\sum_{j=1}^{t}
U_{ij}z_{j}\in R[U_{ij}]$, where $U_{ij}$ are variables. If $\proj
({\SI})$ has rational singularities then $R[U_{ij}]/((g_{1},\ldots
,g_{s}):I)$ has rational singularities. This is in particular the case
when $I$ is generated by a regular sequence and $R/I$ has rational
singularities.} 
\medskip

{\it Proof.} The first statement is a 
combination of Proposition 3.4 and Proposition 3.6(iii). For the
second one, Proposition 3.6 (v) and (vi) imply that ${\RI}\simeq
{\SI}$ and  that ${\cal G}_{I}={\RI}\otimes_{R}R/I\simeq
S\otimes_{R}R/I$ has rational  singularities. Now $\proj ({\cal
G}_{I})$ is a Cartier divisor in  $\proj ({\RI})$, so that the result
is a consequence of the following classical lemma.\fini\medskip 

{\nrm Lemma 3.11.} {\it Let $(R,\im )$ be a local ring
essentially of finite type over a field and $x$ an $R$-regular
element. If $R/xR$ is a ring of rational type, then so is $R$.}\medskip
 
This is proved in [HH1, 4.2(h)] in positive characteristic. In
characteristic zero, it is proved in [Elk, proof of Th{\'e}or{\`e}me 2], and
can also be reduced to the case of positive characteristic.\medskip 

{\it Remark 3.12.} In  Theorem 3.10, if ${\cal G}_{I}$ has rational 
singularities then so does $R/I$ by [Bo]. Also if $z_{1},\ldots
,z_{t}$ is a regular sequence, $\proj ({\cal G}_{I})$ is regular if
and only if $R/I$ is regular.\medskip

In any characteristic, we have:\medskip

{\nrm Theorem 3.13.} {\it Let $R$ be a ring essentially of
finite type over a field, $I:=(z_{1},\ldots ,z_{t})$ an ideal,
and $s\geq 0$ an integer. Let $M$ be a $t$ by $s$ matrix of variables
$U_{ij}$, $S:=R[U_{ij}]$, $(\a_{1}\ldots \a_{s}):=(z_{1}\ldots
z_{t})\cdot M$, $J:=(\a_{1},\ldots ,\a_{s})S:IS$, and $\Ip\in
V(J)$.\smallskip
{\rm (i)} Assume that $\Ip\in V(IS)$, $I_{\Ip \cap R}$ is a complete
intersection, and $(R/I)_{\Ip \cap R}$ is of rational type. Then
$(S/IS+J)_{\Ip}$ is a ring of rational type.\smallskip 
{\rm (ii)} If $\Ip\not\in V(IS)$,  assume that $R_{\Ip\cap R}$ is of
rational type. If $\Ip\in V(IS)$, suppose that $I_{\Ip \cap R}$ is a
complete intersection and $(R/I)_{\Ip \cap R}$ is geometrically of
rational type.  In either case $(S/J)_{\Ip}$ is a ring of rational
type.\smallskip 
}\medskip

{\it Proof.} Recall that (geometric) F-rationality and the rational
singularity property descend and ascend under local 
ring extensions obtained by adjoining finitely many variables and 
localizing ([Ve, 3.1]). This takes care of the case where 
$\Ip\not\in V(IS)$, and we may assume from now on that $\Ip\in V(IS)$.

Localizing at $\Ip \cap R$, we may suppose that $(R,\im ,k)$ is a local
ring with $\im =\Ip \cap R$,  $I=(z_{1},\ldots ,z_{t})$ is a complete 
intersection of height $r$, and $R/I$ is of rational type (in (i)) or 
geometrically of rational type (in (ii)).

After a change of variables over $R$, we may suppose that $z_{i}=0$ 
for $i>r$. We now replace $R$ by $R[U_{ij}\, |\, i>r]_{\Ip\cap 
R[U_{ij}\, |\, i>r]}$. Recall that this does not affect being
(geometrically) of rational type. We may then assume that the
matrix $M$ has size $r$ by $s$. Let $n\geq 0$ be such that $I_{n}(M)
\not\subset \Ip$ and $I_{n+1}(M)\subset \Ip$. Suppose for simplicity 
that the upper left $n$ by $n$ minor of $M$ is not in $\Ip$. 
Again we
replace $R$ by $R[U_{ij}\, |\, i\leq n\; {\rm or}\; j\leq n]_{\Ip\cap
R[U_{ij}\, |\, i\leq n\; {\rm or}\; j\leq n]}$, modify the 
generators of $I$, and perform a change of variables over the ring $R$ 
to assume that $M=\pmatrix{I_{n\times n}&0\cr 0&N\cr}$, where $N$ is
a matrix of variables contained in $\Ip$. Finally, passing to the
ring $R/(z_{1},\ldots ,z_{n})$ we may suppose that $M=N$ is an $r$ 
by $s$ matrix of variables $U_{ij}$ over $R$, $S=R[U_{ij}]$, and
$\Ip =(\im, U_{ij})S$. 

Notice that $J_{\Ip}=(\a_{1},\ldots ,\a_{s},I_{r}(M))S_{\Ip}$ by 
[HU3, 3.4], and that the rings $(S/IS+J)_{\Ip}$ and $(S/J)_{\Ip}$
are domains by [HU3, 3.3]. Moreover, we may assume by induction 
on the dimension that 
$(S/IS+J)_{\Ip}$ and $(S/J)_{\Ip}$ are of rational type locally on the 
punctured spectrum.

As for (i) we consider the natural map 
$$
\phi :R/I\lra B:=(R[U_{ij}]/(I+I_{r}(M)))_{(\im, U_{ij})}=
(S/IS+J)_{\Ip}.
$$
As $R[U_{ij}]/I_{r}(M)$ is $R$-flat, $\phi$ is flat as well. 
Furthermore $\phi$ is local with closed fiber $C:=(k[U_{ij}]/
I_{r}(M))_{(U_{ij})}$.
The ring $C$ is geometrically of rational type by [Ke1, Proposition 2] in 
characteristic 0 and by [HH2, 7.14] in any characteristic. Therefore 
$B$ is of rational type according to [Elk, Th{\'e}or{\`e}me 5] in
characteristic 0 and to [Has, 6.4] or [En, 2.27] in positive characteristic.
To apply the latter, notice that $C$ is geometrically F-injective, 
and that the generic fibre of $\phi$ is F-rational because $B$ is F-rational
on the punctured spectrum and $C$ is F-rational.

To prove (ii) set $A':=k[z_{1},\ldots ,z_{r},U_{ij}]/
(\a_{1},\ldots ,\a_{s},I_{r}(M))$ and let $A$ be the localization of 
$A'$ at the homogeneous maximal ideal. Notice that $z_{1},\ldots ,z_{r}
,U_{ij}$ are indeterminates over $k$. In characteristic 0, $A$ has a 
rational singularity by [Ke2, Ke3]. In positive characteristic, $A'$ is
F-split by [MT, (3), p. 362], hence F-injective. On the other hand $a(A')<0$  
(see for instance [KU2, 2.1]). Thus, since $A$ is F-rational on the 
punctured spectrum, it follows from [FW] (see also [HW, 1.6]) that 
$A$ is F-rational.

Now consider the natural map
$$
\psi : A\lra B:
=(R[U_{ij}]/(\a_{1},\ldots ,\a_{s},I_{r}(M)))_{(\im ,U_{ij})}
=(S/J)_{\Ip}.
$$
Since $R$ is flat over $k[z_{1},\ldots ,z_{r}]_{(z_{1},\ldots
,z_{r})}$, $\psi$ is flat. In addition it is local with 
closed fiber $R/I$, which is geometrically of rational type by
assumption. Again applying  [Elk, Th{\'e}or{\`e}me 5] and [Has, 6.4] or [En,
2.27], we conclude that $B$ is of rational type as well.\fini 

\bigskip

{\bf 4. Castelnuovo-Mumford regularity.}
This last section contains our main result on Castelnuovo-Mumford
regularity. We first prove several facts about the behavior of
local cohomology under liaison. Then we investigate, in the
context of ideals in a standard graded algebra, how several 
properties (such as reducedness, F-rationality, smoothness)
are preserved, or improved, by passing to a generic link
(Theorem 4.4). The main result (Theorem 4.7) will follow from these
facts, combined with Kodaira vanishing theorems and Theorem
1.7.\medskip

{\nrm Proposition 4.1.} {\it Let $k$  be a field, $R$ a positively graded
Gorenstein $k$-algebra with homogeneous maximal ideal $\im$,
$I$, $J$ homogeneous ideals linked by a homogeneous Gorenstein ideal
$\ib $, set $d:=\dim R/\ib$, $a:=a(R/\ib )$, $\Sc :=\proj (R/I)$,
$\Sc' :=\proj (R/J)$ and write $-^{*}$ for the graded
$k$-dual.\smallskip 

{\rm (a)} $H^{d}_{\im}(R/I)\simeq (J/\ib )^{*}[-a]$.\smallskip

{\rm (b)} $H^{d-1}_{\im}(R/I)\simeq {\rm coker}\, (R/J\ra {\rm
End}(\o_{R/J}))^{*}[-a]$. 

Assume $d\geq 2$. Then $ {\rm End}(\o_{R/J})\simeq \oplus_{\mu}H^{0}
(\Sc' _{2},{\cal O}_{\Sc' _{2}}(\mu ))$ where $\Sc' _{2}$ is 
the $S_{2}$-ification of $\Sc' $. Furthermore, $H^{d-1}_{\im}(R/I)$ 
has finite length if and only if $\Sc' $ is $S_{2}$, in which case
$H^{d-1}_{\im}(R/I)\simeq H^{1}_{\im}(R/J)^{*}[-a]$. Finally, if $\Sc' $ 
is reduced over $\overline{k}$, then $\dim_{k}H^{d-1}_{\im}(R/I)_{a}+1$
is the number of components of $\Sc' $ over $\overline{k}$ that are 
connected in codimension $1$.\smallskip

{\rm (c)} $H^{i}_{\im}(R/I)\simeq H^{d-i}_{\im}(R/J)^{*}[-a]$ for
$1\leq i\leq d-1$, if  $\Sc $ satisfies $S_{i+1}$ or $\Sc' $ satisfies
$S_{d-i+1}$.}  
\medskip 
 
{\it Proof.} Part (a) is a consequence of local duality since
$\o_{R/I}\simeq (J/\ib )[a]$. 

To prove (b), we dualize the liaison sequence 
$$
0\ra \omega_{R/J }[-a]\ra  R/\ib \ra R/I \ra 0
$$
into $\o_{R/\ib}$. Using local duality, we obtain an exact sequence
$$
(R/\ib )[a]\lra {\rm End}(\o_{R/J}) [a]\lra \ext^{1}_{R/\ib}(R/I,\o_{R/\ib})
\simeq H^{d-1}_{\im}(R/I)^{*}\lra 0,
$$
which proves the first assertion. As for the last claim notice that
if $k$ is algebraically closed and $\Sc' $ is reduced, then
$\dim_{k}H^{0}(\Sc' _{2},{\cal O}_{\Sc' _{2}})$ is the number 
of connected components of $\Sc' _{2}$, which equals the number of
components of $\Sc' $ connected in codimension $1$.

(c) First assume that $\Sc' $ satisfies $S_{d-i+1}$. For $i=d-1$ our
assertion follows from part (b). Hence we may suppose that $1\leq
i\leq d-2$. Now $\dim \ext^{j}_{R/\ib }(R/J,\o_{R/\ib})\leq i-j-1$ for 
$1\leq j\leq i-1$, and  $\dim \ext^{i}_{R/\ib }(R/J,\o_{R/\ib})\leq 0$. 
Thus we may use the $\o_{R/\ib}$-dual of a free 
$R/\ib$-resolution of $R/J$ to conclude that
$$
H^{i+1}_{\im}(\hom_{R/\ib}(R/J,\o_{R/\ib}))\simeq \ext^{i}_{R/\ib}
(R/J,\o_{R/\ib}).
$$
But the first module is $H^{i+1}_{\im}(\o_{R/J})\simeq H^{i}_{\im}(R/I)[a]$
by the liaison sequence, whereas the latter one is $H^{d-i}_{\im}(R/J)^{*}$.
Finally if $\Sc $ is $S_{i+1}$, our assertion 
follows by reversing the roles of $\Sc $ and  $\Sc' $.\fini\medskip

{\nrm Proposition 4.2.} {\it If the assumptions of Proposition 4.1 are
satisfied then
$$
{\rm reg}(R/I )\leq {\rm reg}(\omega_{R/J}) +a -1.
$$
\indent Moreover, ${\rm reg}(R/I )>d+a-1$ if and only if 
${\rm reg}(\omega_{R/J })\not= d$,
and in this case, ${\rm reg}(R/I )={\rm reg}(\omega_{R/J })+a-1$.
If $d\geq 2$, $\Sc' $ is reduced over $\ol k$ and 
$H^{i}_{\im}(\o_{R/J })_{\geq d-i}=0$ 
for $i\not= d$, then ${\rm reg}(R/I )\leq d+a-2$ if and only if $J_{1}=\ib_{1}$
and $\Sc' $ is connected in codimension 1 over $\overline{k}$.
}\medskip

{\it Proof.} The assertions are easily derived from Proposition 4.1 and 
the liaison sequence (see also [Ch, 5.2]). Also notice that ${\rm
reg}(\omega_{R/J})\geq d$ because $1\in {\rm
End}(\omega_{R/J})_{0}\simeq H^{d}_{\im}(\o_{R/J })_{0}^{*}$.
\fini\medskip 

{\it Definition.} If ${\Sc}$ is an embedded equidimensional projective 
scheme with $\dim \Sc \geq 1$, we
will say that ${\Sc}$ satisfies {\it weak Kodaira vanishing} if $\reg
(\o_{{\Sc}})= \dim \Sc +1$. If ${\Sc}$ is reduced, 
this is equivalent 
to the condition $H^{i}({\Sc},\o_{{\Sc}}(\dim \Sc +1-i))=0$ for 
$1\leq i\leq \dim \Sc -1$.\medskip

{\nrm Corollary 4.3.} {\it Let ${\cal Z}$ be an arithmetically
Gorenstein closed subscheme of a projective space, and let ${\cal S}$
and ${\cal S}'$ be closed subschemes of ${\cal Z}$ of dimension $\geq
1$ linked by forms of degrees $d_{1},\ldots ,d_{r}$. Then ${\rm reg}
({\cal S})<{\rm reg} ({\cal Z})+\sum_{i=1}^{r}(d_{i}-1)$ if and only
if ${\cal S}'$ satisfies weak Kodaira vanishing.} 
\medskip
{\it Proof.} Notice that if ${\cal Z}=\proj (R)$ and $\ib\subset R$ is
the ideal generated by the forms linking $\Sc$ to $\Sc '$, then
$a(R/\ib )=a(R)+d_{1}+\cdots +d_{r}=\reg (R)+\sum_{i=1}^{r}(d_{i}-1)
-\dim (R/\ib )$, 
and therefore Proposition 4.2 gives the assertion.\fini\medskip

{\nrm Theorem 4.4.} {\it Let $k$ be a field, R a standard graded Noetherian 
$k$-algebra, and $I\subset R$ a homogeneous ideal of height $r$ generated
by forms $f_{1},\ldots ,f_{t}$ of degrees $d_{1}\geq \cdots \geq d_{t}
\geq 1$. For an integer $c$ with $r<c<\dim R$, assume that, locally in
codimension $c$, $R$ is Gorenstein and $I$ is a complete
intersection. 

Let $x_{0},\ldots ,x_{n}$ be linear forms spanning $R_{1}$, $a_{ij}
:=\sum_{|\mu |=d_{j}-d_{i}}U_{ij\mu}x^{\mu}$ for $r+1\leq i\leq t$ and
$1\leq j\leq r$, where $U_{ij\mu}$ are variables, write $K:=k(U_{ij\mu})$,
consider the matrix $A:=(a_{ij})$ and define
$$
(\a_{1}\cdots \a_{r}):=(f_{1}\cdots f_{t})\pmatrix{\ {\rm I}_{r\times r}\cr
A\cr}.
$$
Write $R':=R\otimes_{k}K$ and $J:=(\a_{1},\ldots ,\a_{r})R':IR'$.\medskip

{\rm (a)} If $R$ satisfies $S_{r+1}$ and is reduced then $R'/J$ is reduced.
\smallskip

{\rm (b)}  $IR'+J$ has height at least $r+1$ and is a complete
intersection locally in codimension  $c$.\smallskip 

{\rm (c)} $J$ is a complete intersection locally in codimension 
$\min \{ r+3,c\}$.\smallskip

{\rm (d)} If $R/I$ is of rational type locally in codimension $c'<c$,
and locally in codimension $c$ at each prime that does not contain the 
unmixed part of $I$, then $R'/IR'+J$ is of rational type locally in codimension 
$c'+1$, and locally in codimension $c$ at each prime that does not contain 
the unmixed part of $I$. \smallskip

{\rm (e)} Assume that $R/I$ is geometrically of rational type locally 
in codimension $c'<c$. If, locally in codimension $c$ at each prime that 
does not contain the unmixed part of $I$, $R$ is of rational type and $R/I$ 
is geometrically of rational type,  then $R'/J$ is of rational type locally in 
codimension $c'+1$, and locally in codimension 
$c$ at each prime that does not contain the unmixed part of $I$.\smallskip

{\rm (f)} If $R/I$ is smooth (or regular) locally 
in codimension $c'<c$, and $R$ and $R/I$ are smooth (or regular)
locally in codimension $c$ at each prime that 
does not contain the unmixed part of $I$,  then $R'/J$ is smooth (or regular, 
respectively)
locally in codimension $\min \{ r+3,c'+1\}$ at each prime that contains the 
unmixed part of $I$, and locally in codimension 
$\min \{ 2 {\rm ht}\, I_{\Ip\cap R}-r, c\}$ at each prime $\Ip$ that
does not contain the unmixed part of $I$.}\medskip 

{\it Proof.} Let $\Ip \in V(J)$ with $\dim R'_{\Ip}\leq c$ 
and write 
$\ip :=\Ip \cap R$. As $R_{1}\not\subset \ip$, one has $x_{i}\not\in 
\ip$ for some $i$, and then the $a_{ij}$'s are algebraically independent 
over $R_{\ip}$ and $R'_{\ip}$ is a ring of fraction of a polynomial ring 
over $R_{\ip}[a_{ij}]$.
Since the  properties we are interested in are preserved by adjoining
variables and localizing (see the proof of Theorem 3.13), we may
replace $R'$ by $R_{a}:=R_{\ip}[a_{ij}]$ and $\Ip$ by $\Ip_{a}:=\Ip
\cap R_{a}$. Notice that ${\rm ht}\Ip_{a}\leq {\rm ht}\Ip$.  
Let $B=(b_{ij})$ be an $r$ by $r$ matrix 
of variables. Replacing $R_{a}$ by $R_{b}:=R_{a}[b_{ij}]$ and $\Ip_{a}$ 
by $\Ip_{b}:=\Ip_{a}R_{b}$, and multiplying the matrix 
$\pmatrix{\ {\rm I}_{r\times r}\cr A\cr}$ by $B$, 
we may assume that $(\a_{1}\cdots \a_{r}):=(f_{1}\cdots f_{t})\cdot C$,
where $C=\pmatrix{B\cr AB\cr}$ is a matrix of variables over $R_{\ip}$.
Recall that our properties descend under flat local homomorphisms. 
Finally we replace $R'$ by $S:=R_{\ip}[c_{ij}]$, set $\Ip :=\Ip_{b}\cap S$,
and think of $J$ as an ideal of $S$. Notice that 
$R_{\ip}$ is Gorenstein and $I_{\ip}$ is a complete intersection or the 
unit ideal.

As for (a) it suffices to consider the case $\dim S_{\Ip}=c=r$. Now
$I_{\ip}=R_{\ip}$ (see for instance [HU1, the proof of 2.5]), and the
reducedness passes from $R_{\ip}$ to $(S/J)_{\Ip}$.

To prove the other assertions, first assume that ${\rm
ht}I_{\ip}>r$. In this case $J_{\ip}=(\a_{1},\ldots ,\a_{r})_{\ip}$
since $R_{\ip}$ is Cohen-Macaulay. Now (b) and (c) are obvious. For
(d)--(f) we may assume that $I_{\ip}\not=R_{\ip}$. Since $\ip$ cannot
contain the unmixed part of $I$, our assumptions in (d) and (e) imply
that $(R/I)_{\ip}$ is of rational type or geometrically of rational
type, respectively. Now (d) and (e)  follow from Theorem 3.13.  

We now prove (f). After a linear change of variables and deleting 
indeterminates, we may assume that $t=\mu (I_{\ip})={\rm ht} I_{\ip}$. 
Since ${\rm ht} \Ip\leq 2t-r$ and $I\subset \Ip$, it follows that 
$I_{r}(C)\not\subset \Ip$. Thus after inverting a maximal minor of $C$ 
and changing generators of the ideal $IS_{\Ip}$ we may suppose that 
$J_{\Ip}=(f_{1},\ldots ,f_{r})$. Therefore $(S/J)_{\Ip}$ is indeed
smooth (resp. regular).

Finally suppose that ${\rm ht}I_{\ip}=r$. Then $J_{\ip}=L_{1}(I_{\ip})$
is a generic link of $I_{\ip}$. Now (c) follows from [HU1, 2.9(b)]. Also 
${\rm ht}(IS+J)_{\Ip}\geq r+1$ ([HU1, 2.5]), which gives (b) according
to [PS, 1.6] since $(R/I)_{\ip}$ is Gorenstein. 

If $\dim R_{\ip}=\dim S_{\Ip}$,
then $S_{\Ip}=S_{\ip S}$, and therefore $J_{\Ip}=
L_{1}(I_{\ip})_{\ip S}=L^{1}(I_{\ip})$ is a universal link of a 
complete intersection. But then $J_{\Ip}$ is the unit ideal 
by [HU2, 2.13(f)], yielding a contradiction.
Thus $\dim R_{\ip}<\dim S_{\Ip}$. 

Hence if  $\dim S_{\Ip}\leq c'+1$, 
the assumptions of (d) and (e) give that $(R/I)_{\ip}$ is of rational
type or geometrically of rational type, respectively. Again, an application of  
Theorem 3.13 leads to the conclusion.
Finally as to (f), the ring $(R/I)_{\ip}$ is smooth (resp. regular).
If moreover $\dim (S/J)_{\Ip}\leq r+3$, then $(S/J)_{\Ip}$ is
Gorenstein by part (c). Therefore $(S/IS+J)_{\Ip}$ is a complete
intersection on $(S/J)_{\Ip}$, which reduces us to showing that the
former ring is smooth (resp. regular).  
However, after a 
change of variables we may again assume that $t=r$, and then 
$(S/IS+J)_{\Ip}=S_{\Ip}/(IS_{\Ip}+(\det C)S_{\Ip})$, which is smooth
(resp. regular) as the determinant locus is smooth in codimension $3$.\fini
\medskip

{\it Remark 4.5.} Assume $k$ is an infinite field. In Theorem 4.4
parts  (b) and (c), parts (a), (d) and (e) in characteristic zero,
and part (f) for ``smooth'', one can replace the ``generic'' link $J$
by a ``general'' link $J_{u}$ defined by forms $\a_{1}(u),\ldots
,\a_{r}(u)$, with $u$ in a dense open subset of ${\rm Spec}(k[U_{ij\mu
}])$. This is possible due to [Jo, 1.6 and 4.10] and the following
lemma:\medskip

{\nrm Lemma 4.6.} {\it Let $k$ be a field of characteristic zero and
$f:Y\ra U$ a morphism of reduced $k$-schemes of finite type. If
$U$ is irreducible and $Y$ has rational singularities, 
then there exists a non empty open subset $\O$ of $U$ such that 
$Y_{u}=f^{-1}(u)$ has rational singularities for $u\in \O$.}\medskip  

{\it Proof.} Let $Y'\buildrel{\pi}\over{\lra} Y$ be a resolution of
singularities of $Y$. There exists a non empty smooth open subscheme $\O$ 
of $U$ over which $f$ is flat, and such that $\pi$ is a 
simultaneous resolution of singularities over $\O$. If $u\in \O$, 
[Elk, Th{\'e}or{\`e}me 3] implies that $Y_{u}$ has rational 
singularities.\fini\medskip

If $B$ is a standard graded Noetherian algebra over a field, we will
set $B^{top}:=B/\ia$, where 
$\ia$ is the intersection of primary components of $0$ of maximal dimension,
and $B^{low}:=B/H^{0}_{\ia}(B)$. If ${\Sc}:=\proj (B)$, we will set 
accordingly ${\Sc}^{top}:=\proj (B^{top})$ and ${\Sc}^{low}:=\proj (B^{low})$.
Notice that ${\Sc}^{low}=\emptyset$ does not imply ${\Sc}={\Sc}^{top}$. 
Recall that $\irr ({\Sc})$ denotes the locus where ${\Sc}$ is not of rational type.  
\medskip

{\nrm Theorem 4.7.} {\it Let $k$ be a field and let $R$ be a standard graded
Gorenstein $k$-algebra. Let $f_{1},\ldots ,f_{t}$ be forms in $R$ of degrees 
$d_{1}\geq \cdots \geq d_{t}\geq 1$ and $I$ the ideal they generate. 
Set ${\Sc}:=\proj (R/I)$, ${\cal Z}:=\proj (R)$ and $r:={\rm codim}_{{\cal Z}}
{\Sc}$. Suppose $r>0$ and $I$ is not a complete intersection.

\smallskip
{\rm (a)} Assume that one of the following four conditions holds,\smallskip
{\rm (i)} $\dim \Sc =0$,\smallskip 
{\rm (ii)} $\dim \Sc =1$, $\Sc^{top}$ is generically a complete
intersection and  ${\cal Z}$ is 
reduced,\smallskip 
{\rm (iii)} $\dim {\Sc}\in \{ 2,3\}$, ${\rm char}(k)\geq \dim \Sc$,
${\Sc}$ is locally a complete intersection, ${\Sc}^{top}$ is smooth
outside finitely many points and lifts to the ring of second Witt
vectors $W_{2}(k)$,  ${\cal Z}$ and ${\Sc}^{low}$ are smooth  and
$\dim {\Sc}^{low}\leq1$,\smallskip
{\rm (iv)} ${\rm char}(k)=0$, $\dim \irr ({\Sc}^{top})\leq 1$ and, outside 
finitely many points, ${\Sc}$ is locally a complete intersection in ${\cal Z}$ 
and ${\cal Z}$ and ${\Sc}^{low}$ have rational singularities.
\smallskip
Then
$$
\reg ({\Sc}^{top})\leq \reg (R)+d_{1}+\cdots +d_{r}-r-1.
$$
Moreover, if $d_{r}\geq 2$ and $\dim \Sc \geq 1$, assume that $\Sc ''$ is 
a geometrically reduced link of ${\Sc}^{top}$ in some complete intersection
defined by forms of degrees $d_{1},\ldots ,d_{r}$.
Then the inequality  is strict if and only if ${\Sc}''$ is
non-degenerate and geometrically  connected in codimension
one.\smallskip 

{\rm (b)} Assume $d_{1}\geq 2$, $\dim \Sc \geq 1$, $\Sc$ is locally a complete 
intersection in ${\cal Z}$
and one of the following two conditions holds,\smallskip

{\rm (i)} $\Sc$ is of rational type,
\smallskip
{\rm (ii)} ${\rm char}(k)=0$, and 
${\Sc}$ has at most isolated irrational singularities.
\smallskip
Then
$$
\reg (R/I)\leq {{(\dim \Sc +2)!}\over{2}}(\reg (R)+d_{1}+\cdots 
+d_{r}-r-1),
$$
unless $R$ is a polynomial ring over $k$, $I=\ell \, K$ with $\ell$ a 
linear form and $K$ a homogeneous ideal, and either $\dim R=3$ and 
$K$ is a complete intersection of 3 forms of degree $d_1-1$, or else $K$
is generated by linear forms.
}
\medskip

{\it Proof.} Let $J\subset R'$ and $\a_{1},\ldots ,\a_{r}$ be as
in Theorem 4.4, and write $\Sc ':=\proj (R'/J)$.\smallskip

(a) By Proposition 4.2 it suffices to show that $H^{i}_{\im}
(\o_{R'/J})_{>0}=0$ for every $i$.\smallskip
{\rm (i)} The assertion is clear as $R'/J$ is Cohen-Macaulay.\smallskip 
{\rm (ii)} Since $R'/J$ is reduced by Theorem 4.4(a), the claim
follows from the proof of Theorem 1.7(i).\smallskip
{\rm (iii)} Applying Theorem 4.4(f) with $c:=\dim R-1$ and $c':=\dim R-2$, 
we obtain that $\Sc '$ is smooth. Furthermore $\Sc '$ is 
liftable to $W_{2}(k)$, see [Bu, 6.3.10]. Thus indeed [DI, 2.8] 
implies the asserted vanishing.\smallskip
{\rm (iv)} Applying Theorem 4.4(e) with $c:=\dim R-2$ and $c':=\dim R-3$, 
we deduce that $\dim \irr (\Sc ')\leq 0$. Now Theorem 1.3 yields the 
vanishing as claimed.\smallskip

(b) Set $\s :=\reg (R)+\sum_{i=1}^{r}(d_{i}-1)$. First assume that 
$R$ is not regular or $r\geq 2$. In this case $\s \geq d_{r+1}$ and
also $\s \geq 2$, since $d_{1}\geq 2$ and $I$ is not a complete 
intersection. 
Inducting on $\dim {\Sc}$, we may assume that $\dim {\Sc}\geq 2$ 
according to Proposition 2.2. From Theorem 4.4(b) with 
$c:=\dim R-1$ it follows that $R'/IR'+J$ defines a
locally complete intersection scheme $\Y $, necessarily of dimension
$\dim {\Sc}-1$. Also, Theorem 4.4(d) with 
$c:=\dim R-1$ in case (i), $c:=\dim R-2$ in case (ii) and $c':=c-1$
implies that  ${\cal Y}$ is of rational type in case (i) or has at most
isolated irrational singularities in case (ii).
Furthermore we have $\Y =\proj (R'/IR'+(J)_{\leq \s})$ by Theorem
1.7(iii) in the setting of (i) and by Theorem 1.7(ii) in (ii).  
Since $\Y $ is a Cartier divisor in ${\Sc}$ there exist $d:=\dim \Sc +1$
forms $\b_{i}\in J$ of degrees at most $\s$ such that $\Y =\proj (R'/(I,
\b_{1},\ldots ,\b_{d}))$. Set $J':=(\a_{1},\ldots ,\a_{r},\b_{1},
\ldots ,\b_{d})$. One has an exact 
sequence,
$$
0\ra R'/IR'\cap J'\ra R'/IR'\oplus R'/J' \ra R'/IR'+J'\ra 0.
$$
Notice that 
$$
\reg (R/I)=\reg (R'/IR')\leq \max \{ \reg (R'/IR'\cap J'),\reg 
(R'/IR'+J')\} .
$$ 
Now, $IR'\cap J'=(\a_{1},\ldots ,\a_{r})$ so that $\reg (R'/IR'\cap J)
=\s$. Recall that $\s \geq d_{r+1}$. Hence if $IR'+J'$ is not a
complete intersection then our induction hypothesis yields 
$$
\reg (R'/IR'+J')\leq {{d!}\over{2}}(\s -1 +d(\s -1))= {{(d+1)!}\over{2}}
(\s -1).
$$
If on the other hand $IR'+J'$ is a complete intersection then $\reg
(R'/IR'+J')\leq \s +d(\s -1)$. In either case our assertion follows. 

Next, assume $R$ is regular and $r=1$. By Proposition 2.2 we may assume $d
\geq 3$ or $\s =1$. Now $I=\ell \, K$ where $\ell$ is a form of degree
$s\geq 1$ and $K$ a proper homogeneous ideal of dimension at most
one. From Proposition 2.1, applied to $K$, we obtain $\reg (R/I) \leq
s+(d+1)(\s -s)$. But the latter is at most ${{(d+1)!}\over{2}}(\s -1)$
unless $s=\s =1$, in which case $\ell$ is a linear form and $K$ is
generated by linear forms.\fini\medskip 

Even with the assumptions of Theorem 4.7(b), the estimate of Theorem
4.7(a) does no longer hold if one replaces ${\Sc}^{top}$ by $R/I$, as
one can see by taking $R:=k[X,Y,Z]$ and $I:=(X^{3},XY,XZ)$.

\bigskip\bigskip\bigskip

\noindent {\nrm Institut de Math{\'e}matiques, CNRS \&\ Universit{\'e} Paris 6,
4, place Jussieu, \par
\noindent F--75252 Paris {\nrm cedex} 05, France}\par\noindent
{\it E-mail:} chardin@math.jussieu.fr
\bigskip
\noindent {\nrm  Department of Mathematics, Purdue University,
West Lafayette, IN  47907, USA}\par\noindent
{\it E-mail:} ulrich@math.purdue.edu

\vfill\eject

\centerline{\nrm REFERENCES}\bigskip

{\prm
{\prm [BEL] A. Bertram, L. Ein and R. Lazarsfeld: Vanishing theorems, a theorem
of Severi, and the equations defining projective varieties.} {\psl J. Amer. 
Math. Soc.} {\pbo 4} {\prm (1991), 587--602}.\par  

{\prm [Bo] J.-F. Boutot: Singularit{\'e}s rationnelles et quotients par les
groupes r{\'e}ductifs, }{\psl Invent. Math.} {\pbo 88} {\prm (1987), 65--68}.\par

{\prm [BKM] W. Bruns, A. Kustin and M. Miller: The resolution of the generic
residual intersection of a complete intersection, }{\psl J. Algebra} {\pbo
128} {\prm (1990), 214--239}.\par

{\prm [BE] D. Buchsbaum and D. Eisenbud: Remarks on ideals and resolutions,
}{\psl Symp. Math.} {\pbo 11} {\prm (1973), 193--204}.\par

{\prm [Bu] R.-O. Buchweitz: {\pit Contributions {\`a} la th{\'e}orie des
singularit{\'e}s,}  Th{\`e}se de l'Universit{\'e} Paris VII (1981)}.\par

{\prm [Ch] M. Chardin: Applications of some properties of the canonical 
module in computational projective algebraic geometry, }{\psl J. 
Symbolic Comput.} {\pbo 29} {\prm (2000), 527--544}.\par

{\prm [DI] P. Deligne and L. Illusie: Rel{\`e}vements modulo $p^{2}$ et 
d{\'e}composition du complexe de de Rham, }{\psl Invent. Math.} {\pbo 89} 
{\prm (1987), 247--270}.\par

{\prm [Elk] R. Elkik: Singularit{\'e}s rationelles et d{\'e}formations,
}{\psl Invent. Math.} {\pbo 47} {\prm (1978), 139--147}.\par

{\prm [Elz] F. Elzein: Complexe dualisant et applications {\`a} la classe
fondamentale d'un cycle, }{\psl Bull. Soc. Math. France}{\prm , M{\'e}moire}
{\pbo 58} {\prm (1978)}.\par

{\prm [En] F. Enescu: On the behavior of F-rational rings under flat
base change,} {\psl J. Algebra} {\pbo  233} {\prm (2000),  543--566}.\par

{\prm [FW] R. Fedder and K. Watanabe: A characterization of F-regularity in 
terms of F-purity,} {\pit Commutative Algebra  (Berkeley, 1987)}, {\prm
Math. Sci. Res. Inst. Publ. {\pbo 15}, Springer, New York, 1989, 
227--245.}\par

{\prm [Fu] W. Fulton:  {\pit Intersection Theory}, Springer, Berlin, 
1984.}\par

{\prm [Har] N. Hara: A characterization of rational singularities in terms of
injectivity of Frobenius maps, }{\psl Amer. J. Math.} {\pbo 120} {\prm (1998),
981--996}.\par

{\prm [HW] N. Hara and K. Watanabe: The injectivity of Frobenius acting
on cohomology and local cohomology modules, }{\psl Manuscripta Math.}
{\pbo 90} {\prm (1996), 301--315}.\par

{\prm [Has] M. Hashimoto:  Cohen-Macaulay 
F-injective homomorphisms,} {\pit Geometric and Combinatorial Aspects of
Commutative Algebra (Messina, 1999),} {\prm Lecture Notes in Pure and
Appl. Math. {\pbo 217}, Marcel Dekker, 2001, 231--244.}\par 

{\prm [HSV] J. Herzog, A. Simis and W. Vasconcelos: Koszul homology and 
blowing-up rings,} {\pit Commutative Algebra (Proc. Trento 1981), }
{\prm Lecture Notes in Pure and Applied Math. 
{\pbo 84},  Marcel Dekker, New York, 1983, 79--169.}\par

{\prm [Hi] H. Hironaka: Smoothing of algebraic cycles of small dimensions,
}{\psl Amer. J. Math.} {\pbo 90} {\prm (1968), 1--54}.\par 

{\prm [HH1] M. Hochster and C. Huneke:
F-regularity, test elements and smooth base change, }{\psl Trans. Amer.
Math. Soc.} {\pbo 346} {\prm (1994), 1--62}.\par

{\prm [HH2] M. Hochster and C. Huneke: Tight closure of parameter ideals and 
splitting in module-finite extensions, }{\psl J. Algebraic Geom.} {\pbo 3}
{\prm (1994), 599--670}.\par

{\prm [Hu] C. Huneke: On the symmetric and Rees algebras of an ideal 
generated by a $d$-sequence, }{\psl J. Algebra} {\pbo 62} {\prm (1980), 
268--275}.\par 

{\prm [HR] C. Huneke and M. E. Rossi: The dimension and components of 
symmetric algebras, }{\psl J. Algebra} {\pbo 98} {\prm (1986), 200--210}.\par

{\prm [HU1] C. Huneke and B. Ulrich: Divisor class groups and deformations,
}{\psl Amer. J. Math.} {\pbo 107} {\prm (1985), 1265--1303}.\par

{\prm [HU2] C. Huneke and B. Ulrich: The structure of linkage, }{\psl Ann. of 
Math.} {\pbo 126} {\prm (1987), 277--334}.\par

{\prm [HU3] C. Huneke and B. Ulrich: Residual intersections, }{\psl J. reine
angew. Math.} {\pbo 390} {\prm (1988), 1--20}.\par

{\prm [Jo] J.-P. Jouanolou: {\pit Th{\'e}or{\`e}mes de Bertini et applications},
Progress in Math.} {\pbo 42}, {\prm Birkh{\"a}user, Basel, 1983}.\par

{\prm [Ke1] G. Kempf: On the geometry of a theorem of Riemann,
}{\psl Ann. of Math.} {\pbo 98} {\prm (1973), 178--185}.\par 

{\prm [Ke2] G. Kempf: Images of homogeneous vector bundles and
varieties of complexes, }{\psl Bull. Amer. Math. Soc.} {\pbo 81} {\prm
(1975), 900--901}.\par 

{\prm [Ke3] G. Kempf: On the collapsing of homogeneous bundles, 
}{\psl Invent. Math.} {\pbo 37} {\prm (1976), 229--239}.\par

{\prm [Ko] J. Koll{\'a}r: Higher direct images of dualizing sheaves I, }{\psl
Ann. of Math.} {\pbo 123} {\prm (1986), 11--42}.\par

{\prm [KU1] A. Kustin and B. Ulrich: A family of complexes associated to
an almost alternating map, with applications to residual intersections,
}{\psl Mem. Amer. Math. Soc.} {\pbo 461} {\prm (1992)}.\par

{\prm [KU2] A. Kustin and B. Ulrich: If the socle fits,
}{\psl J. Algebra} {\pbo 147} {\prm (1992), 63--80}.\par

{\prm [KW] E. Kunz and R. Waldi: Regular differential forms,} {\pit
Contemporary Math.} {\pbo 79} {\prm (1988)}.\par

{\prm [Li] J. Lipman: Dualizing sheaves, differentials and residues on
algebraic varieties, }{\psl Ast{\'e}risque} {\pbo 117} {\prm (1984)}.\par

{\prm [MS] V. B. Mehta and V. Srinivas: A characterization of rational 
singularities, }{\psl Asian J. Math.} {\pbo 1} {\prm (1997), 249--271}.\par

{\prm [MT] V. B. Mehta and V. Trivedi: Variety of complexes and 
F-splitting, }{\psl J. Algebra} {\pbo 215} {\prm (1999), 352--365}.\par

{\prm [Oh] T. Ohsawa: Vanishing theorems on complete K{\"a}hler manifolds,
}{\psl Publ. Res. Inst. Math. Sci.} {\pbo 20} {\prm (1984), 21--38}.\par

{\prm [PS]  C. Peskine and L. Szpiro: Liaison des vari{\'e}t{\'e}s alg{\'e}briques. 
I, }{\psl Invent. Math.} {\pbo 26} {\prm (1974), 271--302}.\par

{\prm [S1] K. Smith: F-rational rings have rational singularities, }{\psl
Amer. J. Math.} {\pbo 119} {\prm (1997), 159--180}.\par
 
{\prm [S2] K. Smith: Fujita's freeness conjecture in terms of local
cohomology, }{\psl J. Algebraic Geom.} {\pbo 6} {\prm (1997), 417--429}.\par 

{\prm [S3] K. Smith: A tight closure proof of Fujita's freeness conjecture for 
very ample line bundles, }{\psl Math. Ann.} {\pbo 317} {\prm (2000), 
285--293}.\par

{\prm [Va] G. Valla: On the symmetric and Rees algebras of an ideal, 
}{\psl Manuscripta Math.} {\pbo 30} {\prm (1980), 239--255}.\par 

{\prm [Ve] J. Velez: Openness of the F-rational locus and smooth base change,
}{\psl J. Algebra} {\pbo 172} {\prm (1995), 425--453}.\par  

}

\end